\documentclass[11pt]{article}
\usepackage{amsmath, amsthm, amssymb}
\usepackage{geometry}
\usepackage{booktabs}
\usepackage{hyperref}
\usepackage[ruled,vlined]{algorithm2e}
\usepackage{amsmath,amssymb,amsthm}
\usepackage{mathtools}
\usepackage{graphicx}
\usepackage{booktabs}
\usepackage{enumitem}
\usepackage{hyperref}
\usepackage{microtype}
\usepackage{array}
\usepackage{bm}

 \usepackage{tikz}
\usepackage{pgfplots}
\pgfplotsset{compat=1.18}

\newtheorem{theorem}{Theorem}[section]
\newtheorem{proposition}[theorem]{Proposition}
\newtheorem{remark}[theorem]{Remark}
\newtheorem{definition}[theorem]{Definition}
\newtheorem{corollary}[theorem]{Corollary}
\newtheorem{example}[theorem]{Example}

\newcommand{\R}{\mathbb{R}}

\newcommand{\eps}{\varepsilon}

\title{Far-field prediction of algebraic curves via generalized asymptotes}
\author{Sonia P\'erez-D\'{\i}az \\
Dpto. de F\'{\i}sica y Matem\'aticas \\
        Universidad de Alcal\'a \\
      E-28871 Alcal\'a de Henares, Madrid, Spain  \\
sonia.perez@uah.es}
\date{}          

\begin{document}

\maketitle



\begin{abstract}
The asymptotic behavior of algebraic curves is classically described by means of linear asymptotes. However, many infinite branches exhibit a non-linear behavior that cannot be represented accurately by straight lines. In this paper we propose a computational framework for far-field prediction of algebraic curves based on generalized asymptotes. Starting from the asymptotic expansion of a real branch at infinity, we retain its non-decaying terms in order to construct an algebraic surrogate, which is then used as an adaptive far-field model outside a sufficiently large compact region. We introduce explicit discrepancy indicators and a switching radius that determine when the exact branch may be replaced by its generalized asymptotic predictor.  The framework recovers the classical linear asymptote whenever the latter exists, while remaining effective for branches with polynomial or Puiseux-type dominant behavior. In the implicit setting, we prove algebraicity and compatibility properties of the predictor and derive asymptotic error estimates from the first neglected term of the expansion. A residual certificate distinguishes guaranteed finite-tail replacement from asymptotic radius estimates. Tangential tails, slow decay, and a smooth transition rule clarify the limitations of truncation. Reproducible finite-window tests compare the predictors with fitted polynomial baselines; these tests are not a benchmark of a complete plotting system.
\end{abstract}

\\



\section{Introduction}

The study of the behavior of algebraic curves at infinity is a classical topic with direct impact on geometric modeling, qualitative visualization, and symbolic-numeric computation. In practice, the remote part of a curve plays a role in clipping and windowing procedures, adaptive plotting, branch classification, and the construction of simplified surrogate models. A faithful description of that far-field geometry is therefore not only of theoretical interest but also of computational relevance.

The traditional framework for asymptotic analysis is based on linear asymptotes. While this description is useful in many situations, it becomes too restrictive when the dominant remote behavior of a branch is nonlinear. In those cases, replacing a branch by a line may fail to capture the actual geometry at infinity, and therefore provides little information for downstream computational tasks.

This limitation has motivated a series of recent advances on the symbolic computation of asymptotic data for algebraic curves. New algorithms have been proposed for parametric and implicit curves, for asymptotic families, and for generalized asymptotes; see, among others, \cite{SPD2020,SPD2022formula,SPD2022family,SPD2023symmetry,SPD2023amai,SPD2024survey,SPD2024jcam,SPD2025space}. In particular, generalized asymptotes provide a natural extension of the classical notion by allowing the asymptotic model itself to be a low-degree algebraic curve rather than necessarily a straight line.

The computation of Puiseux expansions and the truncation to non-decaying
terms belong to the existing theory cited above. Algebraicity and convergence
of these truncations are recalled here as background, not claimed as new
results. The contribution studied here is a tolerance-based replacement
procedure: a residual test for finite-tail certification, an explicit distinction
between predictor construction and separation of tangential tails, and a smooth
replacement rule. We test its behaviour on specified algebraic examples,
including a deliberately slow-converging family. This does not establish
superiority over all asymptotic or numerical continuation methods.

The scope of this paper is intentionally focused. We do not further develop the general theory of generalized asymptotes. Instead, we exploit it as a practical device for \emph{far-field prediction}. More precisely, for each real branch at infinity of a given algebraic curve, we construct a generalized asymptotic predictor and use it as an adaptive surrogate outside a prescribed compact region. The result is a simple and interpretable model-reduction strategy for the remote geometry of the curve.

The main contributions of the paper are the following:
\begin{enumerate}[label=(\roman*),leftmargin=*]
\item we formulate a far-field prediction framework based on generalized asymptotic surrogates of real infinite branches;
\item we introduce explicit discrepancy indicators and a switching radius that determine when the exact branch can be replaced by its surrogate;
\item we prove uniform far-field convergence and a basic invariance result under affine changes preserving the chosen asymptotic direction;
\item we recall algebraicity and asymptotic error orders, and give a residual-based sufficient certificate for finite-tail replacement;
\item we show that the framework recovers the classical linear asymptote whenever the branch is linearly asymptotic, while extending naturally to genuinely nonlinear dominant behavior;
\item we analyze tangential tails, slow convergence, smooth transitions, and a coordinatewise space-curve extension, and report reproducible finite-window comparisons.
\end{enumerate}

The paper is organized as follows. Section~\ref{sec:predictors} introduces the asymptotic setting, fixes the notation, and develops the notion of generalized asymptotic predictor together with the associated switching criterion. In the same section, Subsection~\ref{sec:linear} discusses the relation with classical linear asymptotes. Section~\ref{sec:puiseux} extends the framework to implicit curves with Puiseux branches, and Subsection~\ref{Sub:PuseuxEx} illustrates this setting with a representative example. Subsection~\ref{sec:examples} collects several exact algebraic examples and explicit discrepancy computations. Section~\ref{sec:application} presents plotting-oriented examples and applications, including a real-world motivation based on the far-field continuation of partially observed algebraic contours. Finally, Section~\ref{sec:conclusions} contains the concluding remarks.

\section{Generalized asymptotic predictors}\label{sec:predictors}

We begin by briefly recalling the asymptotic viewpoint needed in the sequel. Full algorithmic and theoretical details can be found in \cite{SPD2020,SPD2022formula,SPD2022family,SPD2023symmetry,SPD2023amai,SPD2024jcam}.

Let $\mathcal C\subset\R^2$ be a real algebraic curve. A real branch at infinity of $\mathcal C$ is an unbounded connected part of $\mathcal C$ that, after a suitable affine change of coordinates if necessary, admits for large $x$ a representation
\begin{equation}\label{eq:branch}
B=\{(x,\psi(x)):\ x>R_0\},
\end{equation}
where $\psi$ has an asymptotic expansion in powers or Puiseux powers of $x$.

The guiding principle is that the dominant non-decaying terms in that expansion determine a simpler curve with the same far-field behavior.

\begin{definition}
Let $B$ be a real branch at infinity represented as in \eqref{eq:branch}. A \emph{far-field surrogate} of $B$ is a selected real tail $A_B$ of an algebraic curve, represented for sufficiently large $x$ by
\[
A_B=\{(x,\widetilde{\psi}(x)):\ x>R_1\},
\]
and satisfies
\[
\lim_{x\to\infty}\bigl(\psi(x)-\widetilde{\psi}(x)\bigr)=0.
\]
If, moreover, $A_B$ has minimal degree among the algebraic curves with this property, then $A_B$ is called a generalized asymptote of $B$.
\end{definition}

We use the standard non-decaying truncation of the Puiseux expansion as the predictor. We distinguish its selected real tail from the full algebraic curve containing it; minimal-degree characterizations are part of the cited theory and are not reproved here. This operational viewpoint is sufficient for far-field prediction and is consistent with the general framework developed in \cite{SPD2022family,SPD2023symmetry,SPD2023amai}.

\begin{remark}
If the retained terms are affine, namely $\widetilde{\psi}(x)=mx+n$, then the generalized asymptote is the classical linear asymptote. Hence the generalized setting strictly contains the classical one.\\
\end{remark}

Now, let $B$ be a real branch at infinity of a real algebraic curve, represented by
\[
B=\{(x,\psi(x)):\ x>R_0\}.
\]
Assume that $\psi$ admits an asymptotic expansion
\begin{equation}\label{eq:expansion}
\psi(x)=\sum_{k=1}^{m} a_k x^{\alpha_k}+r(x),
\end{equation}
where
\[
\alpha_1>\alpha_2>\cdots>\alpha_m\ge 0,
\]
the coefficients $a_k$ are real, and $r(x)\to 0$ as $x\to\infty$. The exponents may be integers or rational numbers.

\begin{definition}
The \emph{generalized asymptotic predictor} associated with the branch $B$ is the function
\begin{equation}\label{eq:predictor}
\widetilde{\psi}(x)=\sum_{k=1}^{m} a_k x^{\alpha_k}.
\end{equation}
The corresponding algebraic curve
\[
A_B=\{(x,\widetilde{\psi}(x)):\ x>R_1\}
\]
is called the far-field predictor tail of the branch. Its degree means the degree of its algebraic closure, not the degree of a semialgebraic tail.
\end{definition}

This predictor captures exactly the non-decaying part of the asymptotic behavior.

\begin{proposition}\label{prop:convergence}
Let $B$ be a branch at infinity and let $A_B$ be its generalized asymptotic predictor. Then
\[
\lim_{x\to\infty} |\psi(x)-\widetilde{\psi}(x)|=0.
\]
Consequently, for every $\eps>0$ there exists $R_\eps>0$ such that
\[
|\psi(x)-\widetilde{\psi}(x)|<\eps,\qquad x\ge R_\eps.
\]
\end{proposition}

\begin{proof}
By \eqref{eq:expansion} and \eqref{eq:predictor},
\[
\psi(x)-\widetilde{\psi}(x)=r(x),
\]
and therefore the result follows from the fact that $r(x)\to 0$.
\end{proof}

\subsection{Discrepancy indicators}

To use the predictor adaptively, we define the far-field discrepancy indicators
\begin{equation}\label{eq:indicators}
E_\infty(R):=\sup_{x\ge R} |\psi(x)-\widetilde{\psi}(x)|,
\qquad
E_{\rm rel}(R):=\sup_{x\ge R}\frac{|\psi(x)-\widetilde{\psi}(x)|}{1+|\psi(x)|}.
\end{equation}
The first one is an absolute discrepancy, whereas the second provides a scale-aware version.

\begin{proposition}\label{prop:uniform}
Assume that the remainder $r$ in \eqref{eq:expansion} is continuous for $x\ge R_0$ and satisfies $r(x)\to 0$ as $x\to\infty$. Then, for every $\eps>0$, there exists $R_\eps\ge R_0$ such that
\[
\sup_{x\ge R_\eps} |\psi(x)-\widetilde{\psi}(x)|\le \eps.
\]
Equivalently,
\[
\lim_{R\to\infty} E_\infty(R)=0.
\]
\end{proposition}

\begin{proof}
Since $r(x)\to 0$, for every $\eps>0$ there exists $R_\eps$ such that $|r(x)|\le \eps$ for all $x\ge R_\eps$. Using again $\psi(x)-\widetilde{\psi}(x)=r(x)$, one gets
\[
\sup_{x\ge R_\eps} |\psi(x)-\widetilde{\psi}(x)|
=
\sup_{x\ge R_\eps} |r(x)|
\le \eps.
\]
This proves the claim.
\end{proof}

\begin{corollary}\label{cor:monotone}
Under the assumptions of Proposition~\ref{prop:convergence}, if $|\psi(x)-\widetilde{\psi}(x)|$ is eventually monotone decreasing, then
\[
E_\infty(R)=|\psi(R)-\widetilde{\psi}(R)|
\]
for all sufficiently large $R$.
\end{corollary}

\begin{proof}
If the discrepancy is monotone decreasing for $x\ge R$, its supremum on $[R,\infty)$ is attained at $x=R$.
\end{proof}

\begin{proposition}\label{prop:affine}
Let
\[
T(x,y)=(ax+b,cy+d),
\qquad a>0,\ c\neq 0,
\]
be an affine change of variables, with $X=ax+b$, preserving the direction $x\to+\infty$. Let $B$ be a branch represented by $y=\psi(x)$ and let $\widetilde{\psi}$ be its generalized asymptotic predictor. Then the image branch $T(B)$ admits the algebraic surrogate
\[
\widetilde{\psi}_T(X)=c\,\widetilde{\psi}\!\left(\frac{X-b}{a}\right)+d,
\]
and
\[
\lim_{X\to\infty}\Big|\Psi(X)-\widetilde{\psi}_T(X)\Big|=0,
\]
where $\Psi$ denotes the exact branch function of $T(B)$. The transported surrogate need not equal the canonical non-decaying truncation after a translation in $x$; the assertion concerns convergence of the transported tail.
\end{proposition}

\begin{proof}
Writing $X=ax+b$, one has $x=(X-b)/a$ and
\[
Y=c\,\psi\!\left(\frac{X-b}{a}\right)+d.
\]
Hence the transformed exact branch is
\[
\Psi(X)=c\,\psi\!\left(\frac{X-b}{a}\right)+d.
\]
Defining
\[
\widetilde{\psi}_T(X)=c\,\widetilde{\psi}\!\left(\frac{X-b}{a}\right)+d,
\]
we obtain
\[
\Psi(X)-\widetilde{\psi}_T(X)
=
c\left(\psi\!\left(\frac{X-b}{a}\right)-\widetilde{\psi}\!\left(\frac{X-b}{a}\right)\right).
\]
Since $(X-b)/a\to\infty$ as $X\to\infty$, Proposition~\ref{prop:convergence} yields the result.
\end{proof}

\subsection{Switching radius}

Given a tolerance $\eps>0$, we define the \emph{switching radius}
\begin{equation}\label{eq:switch}
R_\eps:=\inf\{R\ge R_0:\ E_\infty(R)\le \eps\}.
\end{equation}
Choose $R_0>0$ on a closed regular tail on which the branch extends continuously. The continuous remainder implies that the infimum is attained whenever it is finite. In computation we use a certified sufficient radius $\widehat R_\eps\ge R_\eps$, rather than assume this infimum can be evaluated exactly.

This leads to a simple computational scheme.

\bigskip
\noindent\textbf{Algorithm 1. Generalized far-field prediction of a real branch}

\medskip
\noindent
\textbf{Input:}
A real algebraic curve $\mathcal C$, a real branch at infinity $B$ of $\mathcal C$, and a prescribed tolerance $\eps>0$.

\smallskip
\noindent
\textbf{Output:}
A generalized asymptotic predictor $A_B$ for the branch $B$, together with a switching radius $R_\eps$ beyond which the predictor approximates the branch within the required tolerance.

\begin{enumerate}[leftmargin=*]
\item Determine a real branch at infinity $B$ of $\mathcal C$ and represent it, after a suitable affine change of coordinates if necessary, in the form
\[
B=\{(x,\psi(x)):\ x>R_0\}.
\]

\item Compute an asymptotic expansion of $\psi(x)$ as $x\to\infty$,
\[
\psi(x)=\sum_{k=1}^{m} a_k x^{\alpha_k}+r(x),
\qquad
\alpha_1>\alpha_2>\cdots>\alpha_m\ge 0,
\]
where $r(x)\to 0$.

\item Construct the generalized asymptotic predictor by retaining the non-decaying terms:
\[
\widetilde{\psi}(x)=\sum_{k=1}^{m} a_k x^{\alpha_k}.
\]
Let $A_B$ denote the algebraic curve associated with $\widetilde{\psi}$.

\item Obtain a proved tail bound $|\psi(x)-\widetilde\psi(x)|\le M(x)$,
valid for every $x\ge R_{\rm cert}$, using an explicit remainder estimate
or Theorem~\ref{thm:tube}. A sampled discrepancy is diagnostic only.
\item Find a sufficient $\widehat R_\eps\ge R_{\rm cert}$ such that
$\sup_{x\ge\widehat R_\eps}M(x)\le\eps$. If verification fails, keep the
exact branch and report that replacement has not been certified.
\item Replace the exact branch $B$ by its generalized asymptotic predictor $A_B$ for all $x\ge\widehat R_\eps$.
That is, use the exact branch in the near field and the predictor in the far field.
\end{enumerate}

The symbolic core of the method is the computation of branches and asymptotic expansions, for which effective algorithms are already available in the literature \cite{SPD2020,SPD2022formula,SPD2022family,SPD2023symmetry,SPD2023amai,SPD2024jcam}.

\begin{example}
Consider the implicit algebraic curve
\[
\mathcal C:\ y^3-xy+1=0.
\]
For its positive real branch at infinity, a dominant-balance argument gives the asymptotic expansion
\[
\psi(x)=x^{1/2}-\frac{1}{2}x^{-1}+O(x^{-5/2}).
\]
Hence the generalized asymptotic predictor is
\[
\widetilde{\psi}(x)=x^{1/2},
\]
which corresponds to the branch $y>0$ of the parabola
\[
y^2=x.
\]
Moreover,
\[
\psi(x)-\widetilde{\psi}(x)\sim -\frac{1}{2}x^{-1},
\qquad
E_\infty(R)\sim \frac{1}{2R}.
\]
If $\eps=0.1$, then the switching radius is approximately
\[
R_\eps\approx 5.
\]
This is only an initial estimate. Example~\ref{ex:cubiccert} shows that $5$ does not meet the tolerance and proves that $6$ is a sufficient switching radius.
\end{example}

\begin{figure}[htbp]
\centering
\begin{tikzpicture}
\begin{axis}[
    width=0.82\textwidth,
    height=0.48\textwidth,
    xlabel={$x$},
    ylabel={$y$},
    xmin=0, xmax=10,
    ymin=0, ymax=3.2,
    samples=250,
    domain=0.8:3.1,
    legend style={at={(0.02,0.98)},anchor=north west},
    grid=major
]

\addplot[very thick] ({x^2 + 1/x}, {x});
\addlegendentry{Exact branch of $y^3-xy+1=0$}

\addplot[dashed, very thick] ({x^2}, {x});
\addlegendentry{Predictor $y=\sqrt{x}$}

\addplot[dotted, thick] coordinates {(6,0) (6,3.2)};
\addlegendentry{Certified switching line $x=6$}

\end{axis}
\end{tikzpicture}
\caption{Illustration of Algorithm 1 for the positive branch of the implicit curve $y^3-xy+1=0$. The exact branch is represented parametrically by $(t^2+t^{-1},t)$, while the generalized asymptotic predictor is the branch of the parabola $y^2=x$ given by $(t^2,t)$. Beyond the switching region, both curves become visually close.}
\label{fig:algorithm_example_implicit}
\end{figure}
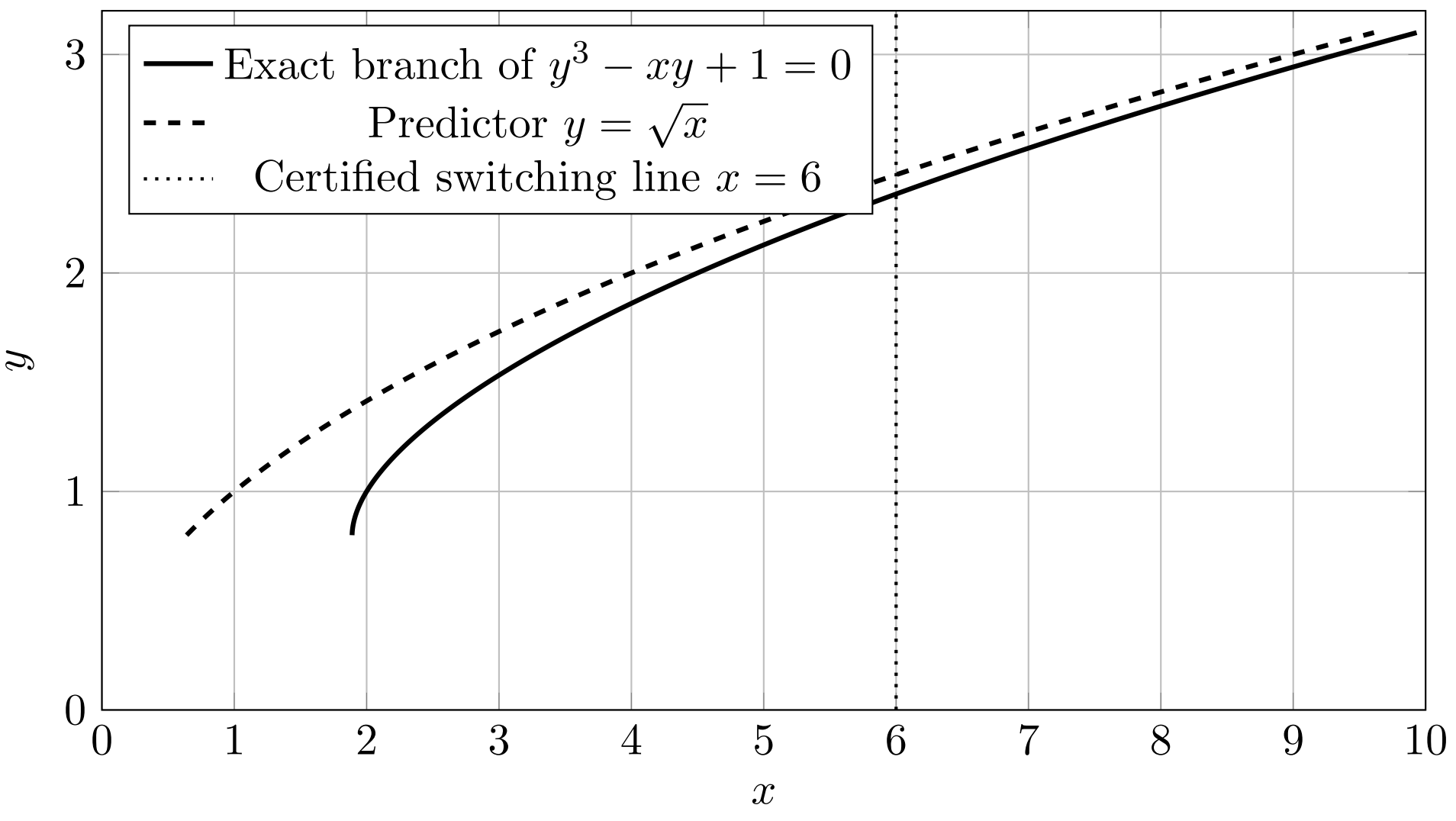

Figure~\ref{fig:algorithm_example_implicit} illustrates the output of Algorithm~1. The exact branch is not introduced through a closed-form expression $y=f(x)$, but through its implicit equation and a parametric representation derived from it. The generalized asymptotic predictor still provides a simple algebraic surrogate, namely the branch of the parabola $y^2=x$, which accurately captures the far-field geometry.

This example also anticipates the more structural implicit Puiseux analysis developed in Section~\ref{sec:puiseux}.

\subsection{Relation with classical linear asymptotes}\label{sec:linear}

The generalized predictor contains the classical theory as a special case.

\begin{proposition}\label{prop:linear}
Suppose that a branch at infinity has an asymptotic expansion
\[
\psi(x)=mx+n+r(x),
\qquad r(x)\to 0 \text{ as } x\to\infty.
\]
Then its generalized asymptotic predictor is the line $y=mx+n$. In particular, the proposed framework recovers the classical linear asymptote whenever it exists.
\end{proposition}

\begin{proof}
The only non-decaying terms in the expansion are $mx+n$, hence the predictor coincides with that affine function.
\end{proof}

The real gain of the generalized framework appears when the dominant behavior is nonlinear.

\begin{proposition}\label{prop:nonlinear}
Suppose that a branch at infinity admits an expansion
\[
\psi(x)=p(x)+r(x),
\]
where $p(x)$ is a nonconstant polynomial of degree at least two and $r(x)\to 0$ as $x\to\infty$. Then no affine function is a valid asymptotic surrogate of the branch, whereas the polynomial curve $y=p(x)$ is a generalized asymptotic predictor.
\end{proposition}

\begin{proof}
Let $\ell(x)=mx+n$ be an affine function. Since $\deg p\ge 2$, one has
\[
|p(x)-\ell(x)|\to\infty
\qquad \text{as }x\to\infty.
\]
Therefore
\[
|\psi(x)-\ell(x)|\ge |p(x)-\ell(x)|-|r(x)|\to\infty,
\]
so $\ell$ cannot define a far-field surrogate. On the other hand,
\[
\psi(x)-p(x)=r(x)\to 0,
\]
hence $y=p(x)$ is a valid predictor.
\end{proof}

Propositions~\ref{prop:linear} and \ref{prop:nonlinear} show that generalized asymptotes are not simply a refinement of linear asymptotes. They enlarge the class of branches that admit a compact asymptotic geometric description, which is precisely what is required in model-reduction tasks at infinity.

\subsection{Finite-tail certification}\label{sec:certificate}
An asymptotic equivalence estimates a switching radius but does not certify a
finite tolerance. The following sufficient test uses the implicit equation.
\begin{theorem}[Residual certificate in an isolating tube]\label{thm:tube}
Let $p$ be a continuous algebraic function on $[R,\infty)$ and $h>0$.
Suppose, for all $x\ge R$, that
\[
 |F(x,p(x))|\le \rho(x),\quad
 |F_y(x,y)|\ge m(x)>0\ (|y-p(x)|\le h),\quad
 \rho(x)<m(x)h.
\]
Then the tube contains exactly one real root $y=\psi(x)$ in each vertical
fibre, and $|\psi(x)-p(x)|\le\rho(x)/m(x)$.
If the selected asymptotic branch is in the tube for sufficiently large $x$,
it is this root branch throughout $[R,\infty)$. Thus
$\sup_{x\ge R}\rho(x)/m(x)\le\eps$ certifies tail replacement.
\end{theorem}
\begin{proof}
At fixed $x$, continuity and nonvanishing of $F_y$ give a constant derivative
sign in the vertical interval. The endpoint values differ from $F(x,p(x))$
in opposite directions by at least $m(x)h$, hence have opposite signs.
The intermediate value theorem gives existence, strict monotonicity gives
uniqueness, and the mean value theorem gives the error bound. The implicit
function theorem supplies a continuous root branch; uniqueness continues its
identification with the selected asymptotic branch across the tail.
\end{proof}
The criterion is sufficient, not necessary. It can fail if the tube contains
several roots or $F_y$ vanishes. A small residual alone does not identify a branch.
For polynomial $F$ with rational or real algebraic coefficients, algebraic $p$,
and semialgebraic $\rho,m$, the inequalities can in principle be verified by
real quantifier elimination, introducing auxiliary algebraic variables and
isolating intervals to specify the real determinations. This finite verification
route may be expensive. Checking a finite sample does not certify an infinite tail.

\begin{example}[A certified radius for the cubic]\label{ex:cubiccert}
For $F=y^3-xy+1$, take $p=\sqrt x$, $R=6$, and $h=1/10$.
On the tube, $F(x,p)=1$ and
\[
 F_y\ge 2x-\tfrac35\sqrt x+\tfrac3{100}.
\]
This lower bound increases for $x\ge6$ and at $6$ exceeds $10$.
Theorem~\ref{thm:tube} certifies error less than $0.1$ on $[6,\infty)$.
The asymptotic estimate $R_{0.1}\approx5$ is not itself sufficient:
$F(5,\sqrt5-0.1)=0.03\sqrt5-0.001>0$ and the large positive root
lies below $\sqrt5-0.1$. Its error at $5$ therefore exceeds $0.1$.
\end{example}

\subsection{From an implicit equation to labelled tails}\label{sec:implementation}
Branch computation uses established Newton--Puiseux methods, not a new
expansion algorithm. The required tasks are as follows.
\begin{enumerate}[leftmargin=*]
\item Remove repeated factors and handle components separately. For a chart
$x\to+\infty$, set $z=1/x$, $y=z^{-k}u$ with integer $k$ large enough
for the relevant branches, and multiply by a power of $z$ to form a polynomial
$H(z,u)$.
\item Compute Newton--Puiseux expansions of $H=0$ at $z=0$.
Each Newton edge gives candidate exponents and a characteristic polynomial
for the leading coefficient. Substitute each root and repeat at multiple
roots until distinct local series separate. Keep determinations real for
$z>0$; remove duplicate parametrizations of the same graph.
Other directions require other charts or $x\mapsto-x$.
\item Return to $y=\psi_i(x)$ and retain every exponent at least zero,
giving $p_i$; if no such term occurs, set $p_i=0$.
Continue farther to separate branch labels or obtain a remainder bound.
These stopping tasks need not have the same truncation depth.
\item At a candidate $R$, isolate the real roots of $F(R,y)$ and match
the selected tail to its isolating interval. Certify a tube as in
Theorem~\ref{thm:tube}, or use an independently proved remainder bound.
Real-root isolation and semialgebraic tests are exact subroutines for
rational or algebraic input.
\item Increase $R$, for example by doubling, until the error certificate
passes. The output is a sufficient $\widehat R_\eps$, not necessarily the
infimum $R_\eps$. If a prescribed search budget is exhausted, report
failure to certify, rather than accept an asymptotic estimate.
\end{enumerate}
This is a verifiable workflow, not a claim of an efficient implementation
for every curve. The numerical tests below implement the stated scalar
examples; they do not implement a general Newton--Puiseux or
quantifier-elimination package.

\subsection{Tangential branches and truncation depth}\label{sec:tangency}
No fixed number of terms suffices for all input curves.
If $\psi_1-\psi_2\sim c x^{-\gamma}$, $\gamma>0$, distinct tails have
the same non-decaying predictor. Their labels require computation through the
first differing coefficient, even though it is not retained in that predictor.
The common predictor alone cannot preserve their separation.

\begin{example}[The parametrization $(t^4+t,t^2)$]\label{ex:tangent}
Elimination gives the irreducible curve
\[
 (x-y^2)^2-y=0.
\]
Its two tails over large positive $x$ satisfy
$x=y^2+\sigma\sqrt y$, $\sigma=\pm1$, $y>0$, and
\[
 \psi_\sigma(x)=x^{1/2}-\frac{\sigma}{2}x^{-1/4}
                         +O(x^{-7/4}).
\]
Indeed, substituting $y=x^{1/2}(1+v)$ balances $2v$ against
$\sigma x^{-3/4}$; the next coefficient in $v$ vanishes.
Both predictors are $p=x^{1/2}$, with curve $y^2=x$.
In the projective chart $X=1$, the point $Y=Z=0$ has common
tangent $Z=0$ for these tails. Computing the $x^{-1/4}$ term
separates them; retaining that decaying term as a corrected surrogate
improves the error order to $x^{-7/4}$. The corrected algebraic
surrogate is not asserted to be a minimal-degree generalized asymptote.
\end{example}

\subsection{Slow decay and smooth transitions}\label{sec:slow}
For a certified error bound $Cx^{-\gamma}$ valid beyond $R_{\rm cert}$,
a sufficient radius is $\max(R_{\rm cert},(C/\eps)^{1/\gamma})$.
Small $\gamma$ may make this impractical. For example,
\[
 F=x(y-1)^{20}-1,\qquad \psi_\pm=1\pm x^{-1/20},\qquad p_\pm=1
\]
give the exact radius $R_\eps=\eps^{-20}$ on $x>0$:
$R_{0.1}=10^{20}$.
Possible responses are to retain additional decaying terms and recertify,
use a relative tolerance when justified, or keep the exact branch within
the window of interest. Here one extra term is exact. In general the cost
and algebraic degree can increase; no universal acceleration is claimed.

A hard switch need not be continuous. Choose $L>0$ and a certified $R$.
For $w(u)=6u^5-15u^4+10u^3$, put
\[
 q(x)=\begin{cases}
 \psi(x),&x\le R,\\
 (1-w((x-R)/L))\psi(x)+w((x-R)/L)p(x),&R<x<R+L,\\
 p(x),&x\ge R+L.
 \end{cases}
\]
Since $0\le w\le1$ and its first two derivatives vanish at both
endpoints, $q$ is $C^2$ if $\psi,p$ are $C^2$ and
$|q-\psi|\le\eps$ on the transition and tail.
Tangent and curvature accuracy require separate derivative bounds.
Exact evaluation remains necessary inside the transition. The composite
model is piecewise algebraic, not generally one algebraic curve.

\subsection{Coordinatewise extension to space curves}\label{sec:spaceextension}
For a given tail $(x,\psi_2(x),\psi_3(x))$ with convergent Puiseux
coordinates, truncate each coordinate to its non-decaying part $p_j$.
A common denominator substitution $x=s^N$ yields a polynomial
parametrization and hence an algebraic space-curve surrogate.
If $|\psi_j-p_j|\le M_j(x)$, its Euclidean error is at most
$\sqrt{M_2(x)^2+M_3(x)^2}$.
For example $(x,x^{1/2}+x^{-1},x^{1/3}+x^{-2})$ has predictor
$(s^6,s^3,s^2)$, $s>0$, and tail error bounded by
$\sqrt{R^{-2}+R^{-4}}$.
This observation does not solve the extraction of space-curve branches from
a general implicit system.

\section{Implicit Puiseux branches and algebraicity of the predictor}\label{sec:puiseux}

The previous sections formulate the far-field predictor assuming that an asymptotic expansion of the branch is already available. We now add a more structural result showing that the framework also applies naturally to \emph{implicit} algebraic curves whose branches at infinity are described by Puiseux expansions.

Consider a real algebraic curve
\[
\mathcal C=\{(x,y)\in\R^2:\ F(x,y)=0\},
\qquad F\in\R[x,y].
\]
Assume that one of its real branches at infinity admits a Puiseux expansion
\begin{equation}\label{eq:puiseux}
\psi(x)=\sum_{j=0}^{\infty} c_j x^{\beta_j},
\qquad
\beta_0>\beta_1>\beta_2>\cdots,
\end{equation}
where the exponents $\beta_j$ are rational numbers and only finitely many of them are nonnegative. Let
\[
\beta_0>\cdots>\beta_m\ge 0>\beta_{m+1}>\beta_{m+2}>\cdots.
\]
We define the corresponding truncation
\begin{equation}\label{eq:puiseuxtrunc}
\widetilde{\psi}(x)=\sum_{j=0}^{m} c_j x^{\beta_j}.
\end{equation}

The next theorem gives two properties that strengthen the mathematical content of the predictor construction: the predictor remains algebraic even in the Puiseux setting, and the first neglected exponent controls the asymptotic order of the error.

\begin{theorem}\label{thm:puiseux}
Let $B$ be a real branch at infinity of an implicit algebraic curve $F(x,y)=0$ admitting a Puiseux expansion \eqref{eq:puiseux}, and let $\widetilde{\psi}$ be the truncation \eqref{eq:puiseuxtrunc}. Then:
\begin{enumerate}[label=(\alph*),leftmargin=*]
\item the graph of $\widetilde{\psi}$ is contained in a real algebraic curve;
\item one has
\[
\psi(x)-\widetilde{\psi}(x)=O(x^{\beta_{m+1}})
\qquad\text{as }x\to\infty,
\]
and therefore
\[
E_\infty(R)=O(R^{\beta_{m+1}})
\qquad\text{as }R\to\infty.
\]
\end{enumerate}
\end{theorem}

\begin{proof}
Let $N$ be a common positive denominator of the rational exponents $\beta_0,\ldots,\beta_m$. Then there exist integers $q_0,\ldots,q_m$ such that $\beta_j=q_j/N$. If we set
\[
x=s^N,
\]
then
\[
\widetilde{\psi}(x)=\sum_{j=0}^{m} c_j s^{q_j}=:P(s).
\]
Hence the predictor is parametrized by
\[
x=s^N,\qquad y=P(s).
\]
Eliminating $s$ from these equations, for instance by taking the resultant of $X-s^N$ and $Y-P(s)$ with respect to $s$, yields a nonzero polynomial relation
\[
G(X,Y)=0
\]
with real coefficients. Therefore the graph of $\widetilde{\psi}$ is contained in a real algebraic curve.

For the error estimate, write
\[
\psi(x)-\widetilde{\psi}(x)=\sum_{j=m+1}^{\infty} c_j x^{\beta_j}.
\]
Under the usual asymptotic interpretation of the Puiseux expansion, the first neglected term determines the order of the remainder. Therefore,
\[
\psi(x)-\widetilde{\psi}(x)=O(x^{\beta_{m+1}}),
\]
which immediately implies
\[
E_\infty(R)=\sup_{x\ge R}|\psi(x)-\widetilde{\psi}(x)|=O(R^{\beta_{m+1}})
\]
as $R\to\infty$.
\end{proof}

Theorem~\ref{thm:puiseux} is relevant for the intended computational use of the predictor. It shows that even when the original branch comes from an implicit curve and is only accessible through local Puiseux analysis, the retained non-decaying terms still define an algebraic far-field surrogate. Moreover, the first discarded exponent provides an immediate asymptotic estimate for the decay of the discrepancy.

\begin{remark}
If the nonnegative exponents in \eqref{eq:puiseux} are integers, then the predictor is polynomial. If fractional nonnegative exponents occur, the predictor remains algebraic but is no longer polynomial; this is precisely the situation in which generalized asymptotes go beyond the polynomial setting while preserving algebraicity.
\end{remark}

The previous theorem can be sharpened in the generic case in which the first discarded coefficient is nonzero.

\begin{corollary}\label{cor:sharp}
Under the assumptions of Theorem~\ref{thm:puiseux}, and assuming that the Puiseux expansion is understood asymptotically in the usual sense, suppose in addition that $c_{m+1}\neq 0$. Then
\[
\psi(x)-\widetilde{\psi}(x)\sim c_{m+1}x^{\beta_{m+1}}
\qquad\text{as }x\to\infty,
\]
and, if the remainder is eventually monotone in absolute value,
\[
E_\infty(R)\sim |c_{m+1}|R^{\beta_{m+1}}
\qquad\text{as }R\to\infty.
\]
\end{corollary}

\begin{proof}
By separating the first neglected term from the rest, one gets
\[
\psi(x)-\widetilde{\psi}(x)
=
c_{m+1}x^{\beta_{m+1}}
+
\sum_{j=m+2}^{\infty} c_jx^{\beta_j}.
\]
Since every exponent $\beta_j$ with $j\ge m+2$ is strictly smaller than $\beta_{m+1}$, the tail is $o(x^{\beta_{m+1}})$. Hence the first asymptotic equivalence follows. The second one is then a direct consequence of Corollary~\ref{cor:monotone}.
\end{proof}
 
A further structural consequence is that the predictor is compatible with the algebraic parametrization naturally associated with the Puiseux branch.

\begin{proposition}\label{prop:paramcompat}
Let $\psi$ be as in \eqref{eq:puiseux}, and let $x=s^N$ be a common denominator substitution for the nonnegative exponents of $\widetilde{\psi}$. Then the far-field predictor is exactly the image of the polynomial truncation in the parameter $s$, namely
\[
x=s^N,\qquad y=P(s),
\]
with
\[
P(s)=\sum_{j=0}^{m} c_js^{N\beta_j}.
\]
In particular, the predictor may be computed either by truncating the branch in $x$ or by truncating its corresponding polynomialized parametrization in $s$.
\end{proposition}

\begin{proof}
This is immediate from the proof of Theorem~\ref{thm:puiseux}. After the substitution $x=s^N$, all retained exponents become integers, so the predictor is represented by the polynomial $P(s)$ in the parameter $s$. Replacing $s$ by $x^{1/N}$ gives back exactly the same truncated function $\widetilde{\psi}(x)$.
\end{proof}

\subsection{An implicit example of Puiseux type}\label{Sub:PuseuxEx}

Consider the irreducible affine curve
\begin{equation}\label{eq:implicitexample}
\mathcal C_5:\ y^3-xy+1=0.
\end{equation}
This curve is naturally given in implicit form, and its far-field behavior is not transparently read off from a global explicit representation. We analyze one of its real branches at infinity.

To determine the dominant balance, we look for a branch of the form $y\sim c x^\alpha$ as $x\to\infty$. Substituting into
\[
y^3-xy+1=0
\]
shows that the dominant terms balance when $3\alpha=\alpha+1$, that is, when $\alpha=1/2$. The corresponding coefficient must satisfy
\[
c^3-c=0,
\]
hence $c=\pm 1$ for the nontrivial branches. Therefore the curve has two nonzero real branches with leading behavior
\[
y\sim \pm x^{1/2}.
\]

We now compute the next term for the positive branch. Write
\[
y=x^{1/2}+u(x),
\qquad u(x)=o(x^{1/2}).
\]
Substituting into \eqref{eq:implicitexample}, we obtain
\[
\left(x^{1/2}+u\right)^3-x\left(x^{1/2}+u\right)+1
=
2xu+3x^{1/2}u^2+u^3+1=0.
\]
The cancellation of the constant term requires
\[
2xu+1=0
\]
at first order, so that
\[
u(x)=-\frac{1}{2x}+O(x^{-5/2}).
\]
Hence the branch admits the Puiseux expansion
\begin{equation}\label{eq:puiseuxexample}
y=x^{1/2}-\frac{1}{2}x^{-1}+O(x^{-5/2}).
\end{equation}

According to the general construction, the far-field predictor is obtained by retaining the nonnegative exponent:
\[
\widetilde{\psi}_5(x)=x^{1/2}.
\]
Therefore the associated generalized asymptote is the algebraic curve
\[
A_5:\ y^2-x=0,
\qquad y>0.
\]
By Theorem~\ref{thm:puiseux}, the predictor is algebraic and
\[
\psi_5(x)-\widetilde{\psi}_5(x)=O(x^{-1}),
\qquad
E_\infty(R)=O(R^{-1}).
\]
In fact, Corollary~\ref{cor:sharp} yields the sharper estimate
\[
\psi_5(x)-\widetilde{\psi}_5(x)\sim -\frac{1}{2}x^{-1},
\qquad
E_\infty(R)\sim \frac{1}{2R}.
\]

This example is mathematically more substantial than the polynomial ones of Subsection~\ref{sec:examples}. The original curve is implicit, the branch is obtained by a Puiseux-type asymptotic argument, and the predictor is an algebraic curve defined only after truncation of the local expansion. Nevertheless, the same far-field philosophy applies unchanged: the remote geometry is captured by the algebraic surrogate
\[
y^2=x.
\]

There is also a third real tail with $\psi_0(x)=x^{-1}+O(x^{-4})$ and predictor $y=0$. The negative outer branch is treated in the same way and leads to the predictor
\[
y=-x^{1/2},
\]
that is, to the lower branch of the same parabola $y^2=x$.

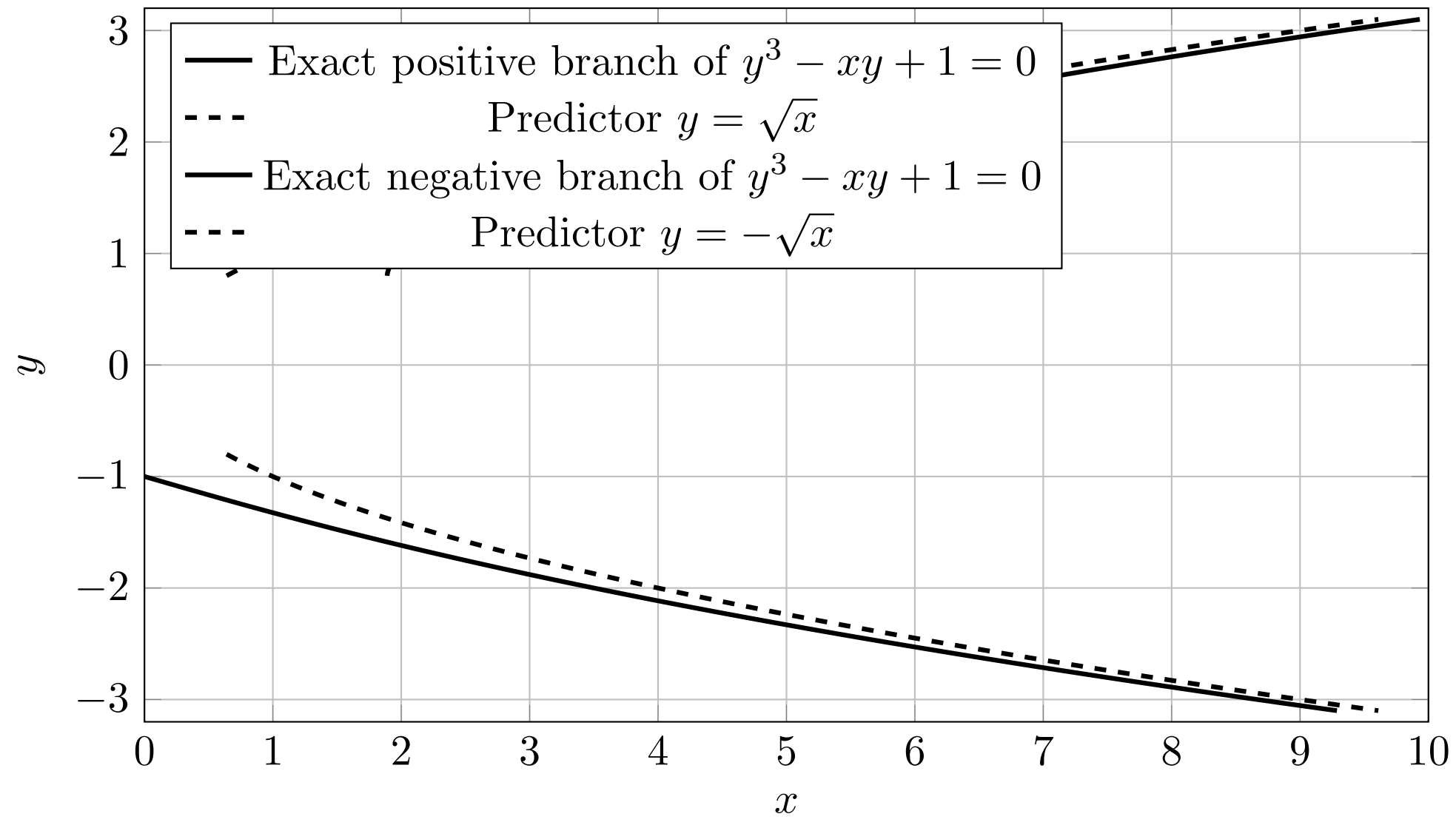
\begin{figure}[htbp]
\centering
\begin{tikzpicture}
\begin{axis}[
    width=0.82\textwidth,
    height=0.50\textwidth,
    xlabel={$x$},
    ylabel={$y$},
    xmin=0, xmax=10,
    ymin=-3.2, ymax=3.2,
    samples=250,
    domain=0.8:3.1,
    legend style={at={(0.02,0.98)},anchor=north west},
    grid=major
]

\addplot[very thick] ({x^2 + 1/x}, {x});
\addlegendentry{Exact positive branch of $y^3-xy+1=0$}

\addplot[dashed, very thick] ({x^2}, {x});
\addlegendentry{Predictor $y=\sqrt{x}$}

\addplot[very thick] ({x^2 - 1/x}, {-x});
\addlegendentry{Exact negative branch of $y^3-xy+1=0$}

\addplot[dashed, very thick] ({x^2}, {-x});
\addlegendentry{Predictor $y=-\sqrt{x}$}

\end{axis}
\end{tikzpicture}
\caption{Far-field behavior of the implicit curve $y^3-xy+1=0$. The outer branches are represented parametrically by $(t^2+t^{-1},t)$ and $(t^2-t^{-1},-t)$, while their generalized asymptotic predictors are the two branches of the parabola $y^2=x$, namely $(t^2,\pm t)$. The coincidence becomes visually clear in the remote region.}
\label{fig:puiseuximplicit}
\end{figure}

Figure~\ref{fig:puiseuximplicit} illustrates the two real nonzero branches of the implicit curve together with their generalized asymptotic predictors. The exact branches are drawn through the parametric representation
\[
(x,y)=(t^2+t^{-1},t)\quad\hbox{or}\quad(t^2-t^{-1},-t),
\]
whereas the asymptotic surrogate is given by the parabola
\[
x=t^2,\qquad y=\pm t.
\]
This visualization makes the far-field approximation transparent and confirms that the algebraic predictor captures the remote geometry of the implicit curve.

\subsection{Exact algebraic examples}\label{sec:examples}

We now present four exact examples. They are intentionally elementary, since the objective is to make the far-field mechanism transparent and to derive closed-form discrepancy expressions.

\begin{example}[{\sf A linear asymptote recovered exactly}] Consider
\[
\mathcal C_1:\ xy-x^2-1=0.
\]
For $x\neq 0$, the curve can be written as
\[
y=x+\frac{1}{x}.
\]
Hence the branch for $x>0$ has predictor
\[
\widetilde{\psi}_1(x)=x,
\]
which is exactly the classical asymptote. The discrepancy is
\[
|\psi_1(x)-\widetilde{\psi}_1(x)|=\frac{1}{x},
\]
and thus
\[
E_\infty(R)=\frac{1}{R},
\qquad
R_\eps=\frac{1}{\eps}.
\]
This example confirms that the generalized framework is consistent with the classical one. \\\end{example}

\begin{example}[{\sf A quadratic far-field surrogate}]

Consider
\[
\mathcal C_2:\ xy-x^3-1=0.
\]
For $x\neq 0$,
\[
y=x^2+\frac{1}{x}.
\]
The predictor is therefore
\[
\widetilde{\psi}_2(x)=x^2.
\]
Moreover,
\[
|\psi_2(x)-\widetilde{\psi}_2(x)|=\frac{1}{x},
\qquad
E_\infty(R)=\frac{1}{R},
\qquad
R_\eps=\frac{1}{\eps}.
\]

The relevant point is not the size of the remainder, but the geometry of the predictor: no line can model the branch at infinity, while the parabola $y=x^2$ captures its dominant behavior exactly.\\
\end{example}

\begin{example}[{\sf Same predictor, different switching radius}]
Consider
\[
\mathcal C_3:\ x^2y-x^4-1=0.
\]
For $x\neq 0$,
\[
y=x^2+\frac{1}{x^2}.
\]
Hence the predictor is again
\[
\widetilde{\psi}_3(x)=x^2,
\]
but now the discrepancy satisfies
\[
|\psi_3(x)-\widetilde{\psi}_3(x)|=\frac{1}{x^2},
\qquad
E_\infty(R)=\frac{1}{R^2},
\qquad
R_\eps=\frac{1}{\sqrt{\eps}}.
\]

This shows that branches with the same generalized asymptote may have different adaptive behavior, and therefore the predictor and the switching radius play complementary roles.\\
\end{example}

\begin{example}[{\sf A cubic predictor}] Consider
\[
\mathcal C_4:\ xy-x^4-1=0.
\]
For $x\neq 0$,
\[
y=x^3+\frac{1}{x}.
\]
Thus
\[
\widetilde{\psi}_4(x)=x^3,
\qquad
E_\infty(R)=\frac{1}{R},
\qquad
R_\eps=\frac{1}{\eps}.
\]
Again, no affine asymptotic model is meaningful, whereas the cubic predictor captures the far-field geometry in a natural way.\\\end{example}
 
Table~\ref{tab:summary} collects the exact predictors and discrepancies obtained above.

\begin{table}[htbp]
\centering
\begin{tabular}{>{\raggedright}p{0.24\textwidth} >{\raggedright}p{0.24\textwidth} >{\centering}p{0.16\textwidth} >{\centering\arraybackslash}p{0.18\textwidth}}
\toprule
Curve & Predictor & $E_\infty(R)$ & $R_\eps$ \\
\midrule
$xy-x^2-1=0$ & $y=x$ & $1/R$ & $1/\eps$ \\
$xy-x^3-1=0$ & $y=x^2$ & $1/R$ & $1/\eps$ \\
$x^2y-x^4-1=0$ & $y=x^2$ & $1/R^2$ & $1/\sqrt{\eps}$ \\
$xy-x^4-1=0$ & $y=x^3$ & $1/R$ & $1/\eps$ \\
\bottomrule
\end{tabular}
\caption{Exact generalized predictors and switching radii for the four examples.}
\label{tab:summary}
\end{table}

The examples make two points clear. First, the method recovers the classical linear asymptote when that model is adequate. Second, in the genuinely nonlinear cases, the generalized predictor yields a compact far-field description that has no linear counterpart.

\section{Some plotting-oriented examples and applications}\label{sec:application}

We next include a simple visual example showing how the generalized asymptotic predictor can be used in adaptive plotting.

Consider again the curve
\[
\mathcal C:\ xy-x^3-1=0,
\]
whose branch for $x>0$ is
\[
y=x^2+\frac{1}{x}.
\]
Its generalized predictor is the parabola
\[
y=x^2.
\]
For this branch, the discrepancy is exactly $1/x$, and hence the switching radius associated with a tolerance $\eps$ is
\[
R_\eps=\frac{1}{\eps}.
\]
For instance, if one prescribes an absolute far-field tolerance $\eps=0.2$, then
\[
R_{0.2}=5.
\]
Therefore, in a plotting procedure, one may evaluate the exact branch on $[x_0,5]$, with $x_0>0$ and replace it by the predictor $y=x^2$ for $x\ge 5$.

This strategy has a direct computational interpretation: the exact branch may involve expensive algebraic manipulations or numerical evaluation steps, while the predictor is a low-degree polynomial that can be sampled essentially at no cost. Thus the surrogate preserves the shape at infinity and, at the same time, simplifies the remote part of the plotting process.

\begin{figure}[htbp]
\centering
\begin{tikzpicture}
\begin{axis}[
    width=0.82\textwidth,
    height=0.48\textwidth,
    xlabel={$x$},
    ylabel={$y$},
    xmin=0.6, xmax=8,
    ymin=0, ymax=66,
    domain=0.6:8,
    samples=200,
    legend style={at={(0.02,0.98)},anchor=north west},
    grid=major
]
\addplot[very thick] {x^2 + 1/x};
\addlegendentry{Exact branch $y=x^2+\frac{1}{x}$}

\addplot[dashed, very thick] {x^2};
\addlegendentry{Generalized predictor $y=x^2$}

\addplot[dotted, thick] coordinates {(5,0) (5,66)};
\addlegendentry{Switching radius $R_{0.2}=5$}
\end{axis}
\end{tikzpicture}
\caption{Far-field prediction for the branch $y=x^2+\frac{1}{x}$. Beyond the switching radius, the generalized asymptote $y=x^2$ provides a visually faithful and computationally simpler surrogate.}
\label{fig:farfield}
\end{figure}
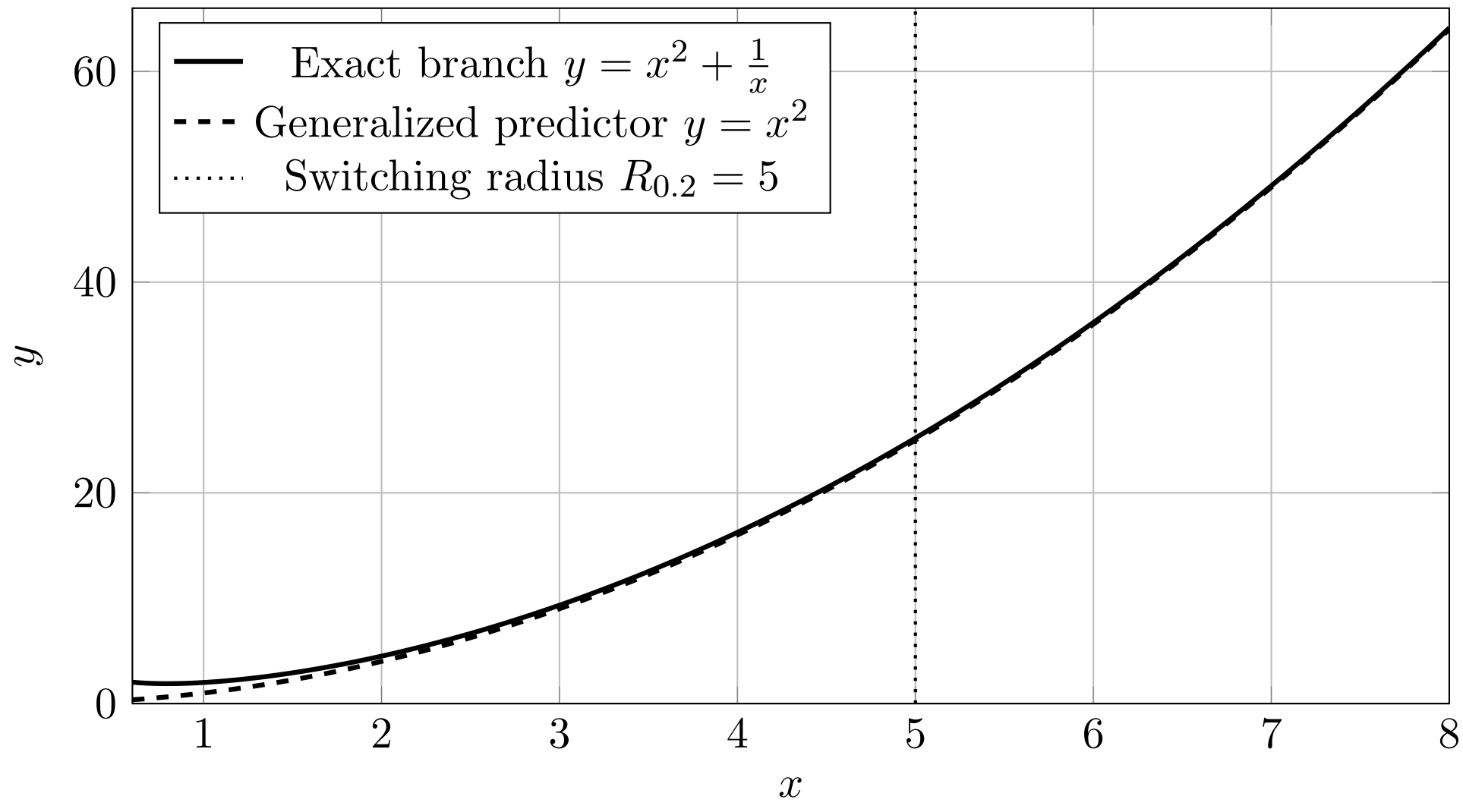

Figure~\ref{fig:farfield} shows the exact branch and its generalized predictor. The two curves are visually almost indistinguishable in the remote region, whereas any affine model would completely miss the correct far-field geometry.

From the viewpoint of adaptive plotting, the rule is therefore clear:
\[
\text{exact branch on } [x_0,R_\eps],\ x_0>0, \qquad
\text{predictor on } [R_\eps,\infty).
\]
This is a minimal but concrete application of generalized asymptotes to geometric processing.

\subsection{An implicit high-degree example with non-polynomial far-field predictor}\label{Sub: predictor}

We next consider an implicit algebraic curve for which the variable $y$ is not conveniently available as an explicit function of $x$, but whose far-field behavior can still be captured by a generalized asymptotic predictor.

Let
\[
\mathcal C_6:\ y^5+y-x^2=0.
\]
This is an implicit algebraic curve. Although, for each sufficiently large $x>0$, the branch we are interested in is locally the graph of a real function, the equation is not naturally handled by solving explicitly for $y$ in terms of $x$. Instead, we extract its far-field behavior directly from the implicit relation.

To determine the dominant balance, we look for a branch of the form
\[
y\sim c\,x^\alpha
\qquad\text{as }x\to\infty.
\]
Substituting into the equation yields
\[
c^5x^{5\alpha}+cx^\alpha-x^2\sim 0.
\]
The dominant terms balance when
\[
5\alpha=2,
\qquad\text{that is,}\qquad
\alpha=\frac{2}{5},
\]
and the leading coefficient must satisfy
\[
c^5=1.
\]
Hence the positive real branch satisfies
\[
y\sim x^{2/5}.
\]

We refine this approximation by writing
\[
y=x^{2/5}(1+u(x)),
\qquad u(x)\to 0.
\]
Substituting into $y^5+y-x^2=0$, we get
\[
x^2\bigl((1+u)^5-1\bigr)+x^{2/5}(1+u)=0.
\]
Since $u$ is small, the first-order balance gives
\[
5x^2u+x^{2/5}\approx 0,
\]
so that
\[
u(x)\approx -\frac{1}{5}x^{-8/5}.
\]
Therefore the branch admits the asymptotic expansion
\[
y=x^{2/5}-\frac{1}{5}x^{-6/5}+O(x^{-14/5}).
\]

According to the general construction, the generalized asymptotic predictor is obtained by retaining the non-decaying term:
\[
\widetilde{\psi}_6(x)=x^{2/5}.
\]
Thus the associated generalized asymptote is the algebraic curve
\[
A_6:\ y^5-x^2=0,
\qquad y>0.
\]
Moreover,
\[
\psi_6(x)-\widetilde{\psi}_6(x)= -\frac{1}{5}x^{-6/5}+O(x^{-14/5}),
\]
and hence
\[
E_\infty(R)=O(R^{-6/5}).
\]

This example is important for two reasons. First, the original curve is implicit and the asymptotic model is extracted without solving it globally for $y$. Second, the far-field predictor is not linear and not polynomial in $x$, but it is still algebraic and therefore perfectly consistent with the generalized asymptotic framework.

For plotting purposes, it is convenient to represent the exact branch parametrically. Setting $y=t$ with $t>0$, the exact branch is given by
\[
x=\sqrt{t^5+t},\qquad y=t,
\]
whereas the predictor is represented by
\[
x=t^{5/2},\qquad y=t.
\]
This allows a direct visual comparison between the exact branch and its far-field surrogate.

As a simple plotting rule, one may choose a tolerance and replace the exact branch by the predictor once the discrepancy becomes small enough. Since the first neglected term has size approximately $(1/5)x^{-6/5}$, the discrepancy decays faster than in the basic examples with $1/x$ decay. The decay exponent alone does not order finite switching radii: constants and the onset of the remainder bound also matter.

\begin{figure}[htbp]
\centering
\begin{tikzpicture}
\begin{axis}[
    width=0.82\textwidth,
    height=0.48\textwidth,
    xlabel={$x$},
    ylabel={$y$},
    xmin=0, xmax=24,
    ymin=0, ymax=3.6,
    samples=250,
    domain=0.2:3.4,
    legend style={at={(0.02,0.98)},anchor=north west},
    grid=major
]
\addplot[very thick] ({sqrt(x^5 + x)}, {x});
\addlegendentry{Exact branch of $y^5+y-x^2=0$}

\addplot[dashed, very thick] ({x^(5/2)}, {x});
\addlegendentry{Generalized predictor $y=x^{2/5}$}

\addplot[dotted, thick] coordinates {(5,0) (5,3.6)};
\addlegendentry{Reference switching line $x=5$}
\end{axis}
\end{tikzpicture}
\caption{Far-field prediction for the positive branch of the implicit curve $y^5+y-x^2=0$. The exact branch is plotted parametrically as $(\sqrt{t^5+t},t)$, while the generalized asymptotic predictor is $(t^{5/2},t)$. In the remote region, both curves become visually indistinguishable.}
\label{fig:implicitfarfield}
\end{figure}
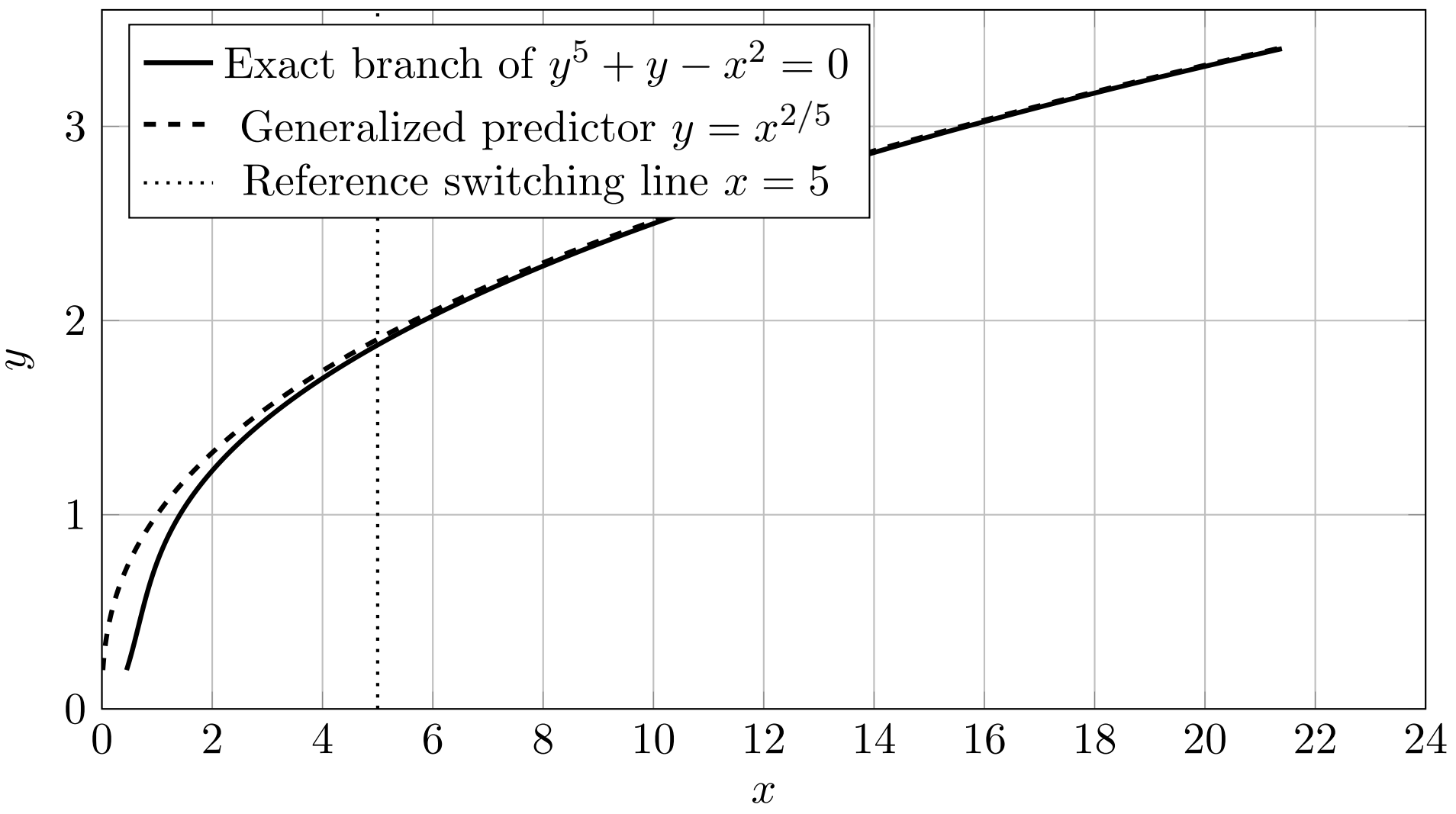

\begin{remark}[{\sf Computational remarks}] The proposed framework is particularly suitable for lightweight computational tasks involving the remote geometry of a curve.\\

\noindent
{\sf Adaptive plotting.}
Once the switching radius has been reached, the exact branch may be replaced by its predictor. This avoids unnecessary evaluation of the original expression while preserving the correct asymptotic shape.\\

\noindent
{\sf Far-field clipping and windowing.}
The predictor provides a simple explicit description of the branch outside a compact set. This can be exploited in clipping procedures and in the construction of better viewing windows for unbounded curves.\\

\noindent
{\sf Qualitative branch analysis.}
The degree and shape of the predictor immediately reveal whether the branch behaves linearly, quadratically, cubically, or according to a more general Puiseux law. This may be useful as a symbolic preprocessing stage before more expensive numerical tasks.\\

These observations support the use of generalized asymptotes as computational surrogates rather than as purely descriptive asymptotic objects.
  
\end{remark}

\subsection{A real-world motivation: far-field continuation of partially observed algebraic contours}\label{sub:realworld}

A natural real-world motivation for generalized far-field prediction arises in geometric vision and contour processing. In many applications, a planar object boundary is modeled, either exactly or approximately, by an implicit algebraic curve
\[
F(x,y)=0.
\]
Typical examples appear in computer-aided inspection, geometric reconstruction from images, and contour-based shape analysis.

In such settings, one often has reliable information only in a bounded observation window. Outside that region, the contour may be partially occluded, truncated by the image frame, or too expensive to track numerically with the full implicit representation. A practical question is then the following: how can one continue the contour in a geometrically meaningful way beyond the observed region?

The framework developed in this paper provides a simple answer. Once a real branch at infinity has been identified, its generalized asymptotic predictor supplies a low-complexity algebraic surrogate that preserves the correct remote geometry of the contour. Therefore, instead of extrapolating the contour by linear heuristics or by purely numerical continuation, one may switch to the generalized asymptotic model beyond a prescribed radius $R_\eps$.

This has two concrete advantages. First, the continuation remains geometrically consistent with the algebraic structure of the original contour. Second, the predictor may be cheaper to evaluate after symbolic preprocessing. No general speed-up is asserted, and for some examples its algebraic degree equals that of the original curve.

For instance, in the implicit example
\[
y^3-xy+1=0,
\]
the visible branch may be represented exactly in the observation region, while its far-field continuation is replaced by the predictor
\[
y^2=x.
\]
In a contour-tracking pipeline, this means that once the branch enters a remote region where only its large-scale behavior matters, the exact model can be replaced by the parabola without losing the essential geometric information.

These are potential uses under the assumption that the implicit model and the selected tail are already known. No recovery of an unknown contour from bounded or noisy observations is proved here, and no image-processing dataset is evaluated. Subject to this restriction, generalized asymptotes may serve as mathematically informed far-field surrogates in:
\begin{enumerate}[label=(\roman*),leftmargin=*]
\item contour continuation outside a finite field of view,
\item continuation of a specified algebraic model beyond a bounded window,
\item geometric preprocessing for clipping, plotting, and visualization of implicit shapes.
\end{enumerate}

We emphasize that the present paper does not address the full image-processing pipeline. Rather, it isolates a geometric core task: the asymptotically faithful continuation of algebraic contours. In this respect, generalized asymptotes provide a mathematically grounded alternative to linear extrapolation in situations where the remote behavior is intrinsically nonlinear.

\subsection{Reproducible finite-window comparisons}\label{sec:validation}
The accompanying dependency-free Node.js script \texttt{validate\_examples.cjs}
generates Table~\ref{tab:validation} and the CSV file.
We use IEEE double arithmetic, $401$ logarithmically spaced abscissae
$x_i=4\,16^{i/400}$ on $[4,64]$, and bisection with $100$ iterations
on explicitly bracketed real roots. Here $d$ is the total degree of the
input curve, $b_+$ the number of real graph tails over $x\to+\infty$
(not the number of projective places), and $d_p$ the degree of the
predictor's algebraic closure.
The labels Tplus/Tminus correspond to $\sigma=+1/-1$ in
Example~\ref{ex:tangent}; Splus/Sminus denote the slow family.
For the cubic, plus/minus denote the two outer roots and small the
root tending to zero.

The measured vertical error is $e_v=\max_i|\psi(x_i)-p(x_i)|$.
We also report the symmetric Hausdorff distance $d_{H,N}$ between the two
finite point sets $\{(x_i,\psi(x_i))\}$ and $\{(x_i,p(x_i))\}$,
computed by exhaustive nearest-neighbour searches in both directions.
It is a sampled-set distance, not a certified Hausdorff distance between
continuous curves or infinite tails. In these tests it equals $e_v$
because of sampling and endpoint effects; it supplies no sharper geometric
certificate. For continuous tails at matching abscissae the actual
Hausdorff distance is at most $E_\infty(R)$, by using the matched point
in each direction.

The baselines are least-squares affine and quadratic polynomials fitted to
$81$ equally spaced exact branch values on $[2,4]$, using the scaled
coordinate $u=x-3$. Their errors $e_{\rm lin}$ and $e_{\rm quad}$ are
measured on the same test window. A fitted affine baseline is not a
classical asymptote; for C1 the true classical asymptote is the predictor.
These baselines measure extrapolation beyond the fitted interval and
do not represent every polynomial fitting strategy. The script specifies
all root brackets and verifies scaled polynomial residuals.
\begin{table}[htbp]
\centering\scriptsize
\begin{tabular}{lrrrrrrr}
\toprule
Tail & $d$ & $b_+$ & $d_p$ & $e_v$ & $d_{H,N}$ & $e_{\rm lin}$ & $e_{\rm quad}$\\
\midrule
C1 & 2 & 1 & 1 & $2.50\cdot10^{-1}$ & $2.50\cdot10^{-1}$ & $6.95\cdot10^{+0}$ & $1.46\cdot10^{+2}$ \\
C2 & 3 & 1 & 2 & $2.50\cdot10^{-1}$ & $2.50\cdot10^{-1}$ & $3.73\cdot10^{+3}$ & $1.46\cdot10^{+2}$ \\
C3 & 4 & 1 & 2 & $6.25\cdot10^{-2}$ & $6.25\cdot10^{-2}$ & $3.73\cdot10^{+3}$ & $1.59\cdot10^{+2}$ \\
C4 & 4 & 1 & 3 & $2.50\cdot10^{-1}$ & $2.50\cdot10^{-1}$ & $2.60\cdot10^{+5}$ & $2.27\cdot10^{+5}$ \\
C5plus & 3 & 3 & 2 & $1.39\cdot10^{-1}$ & $1.39\cdot10^{-1}$ & $1.82\cdot10^{+1}$ & $3.31\cdot10^{+2}$ \\
C5minus & 3 & 3 & 2 & $1.15\cdot10^{-1}$ & $1.15\cdot10^{-1}$ & $9.00\cdot10^{+0}$ & $3.90\cdot10^{+1}$ \\
C5small & 3 & 3 & 1 & $2.54\cdot10^{-1}$ & $2.54\cdot10^{-1}$ & $9.22\cdot10^{+0}$ & $2.92\cdot10^{+2}$ \\
C6 & 5 & 1 & 5 & $3.87\cdot10^{-2}$ & $3.87\cdot10^{-2}$ & $1.05\cdot10^{+1}$ & $1.08\cdot10^{+2}$ \\
Tplus & 4 & 2 & 2 & $3.52\cdot10^{-1}$ & $3.52\cdot10^{-1}$ & $1.32\cdot10^{+1}$ & $1.01\cdot10^{+2}$ \\
Tminus & 4 & 2 & 2 & $3.52\cdot10^{-1}$ & $3.52\cdot10^{-1}$ & $9.76\cdot10^{+0}$ & $5.83\cdot10^{+1}$ \\
Splus & 21 & 2 & 1 & $9.33\cdot10^{-1}$ & $9.33\cdot10^{-1}$ & $8.52\cdot10^{-1}$ & $9.98\cdot10^{+0}$ \\
Sminus & 21 & 2 & 1 & $9.33\cdot10^{-1}$ & $9.33\cdot10^{-1}$ & $8.52\cdot10^{-1}$ & $9.98\cdot10^{+0}$ \\
\bottomrule
\end{tabular}
\caption{Deterministic finite-window tests on $[4,64]$. Errors are sampled measurements, not certified tail errors.}
\label{tab:validation}
\end{table}

The non-decaying predictor is substantially more accurate than these two
fitted baselines for C1--C6 and the tangential example. This reflects the
use of exact far-field information and is not a fair claim of superiority
for recovering unknown curves from observations alone. On the slow family,
the fitted line has a smaller measured error than the uncorrected predictor.
The predictor degrees also show that reduction of algebraic degree is not
universal: C6 has input and predictor degree five. No runtime benefit or
general robustness is inferred from this small deterministic experiment.

\section{Conclusions}\label{sec:conclusions}
Generalized asymptotic truncations provide algebraic models for specified real
tails. Their asymptotic convergence is classical; finite-tolerance replacement
requires additional information. We have separated the asymptotic radius
estimate from a sufficient residual certificate, distinguished branch
separation from predictor truncation, and described smooth replacement
with controlled positional error. Tangential and slow-converging examples
show why these distinctions matter.

The tests concern exact algebraic input and finite windows. They demonstrate
specific error behaviour, not a universal gain in runtime or a validated
vision application. The slow example shows that the minimal truncation can
be impractical and can lose to fitted baselines on a bounded window.
Symbolic preprocessing, real-branch extraction and certification can be
costly. Vertical discrepancies depend on the coordinate chart; they bound
paired-tail Euclidean Hausdorff distance from above but need not approximate
it sharply. A common predictor can merge distinct tangential tails.
The smooth transition controls position, not tangent or curvature errors.
The space-curve observation assumes that a compatible branch expansion is
already available. General branch extraction, runtime benchmarks and
data-driven contour recovery remain outside the scope of this paper.

\end{document}